\documentclass[12pt]{article}
\usepackage[hmargin=.75in,vmargin=1in]{geometry}
\usepackage[american]{babel}
\usepackage[T1]{fontenc}
\usepackage{times}
\usepackage{caption}
\usepackage{textcomp}
\usepackage{epsfig,graphicx}
\usepackage{xcolor}
\usepackage{amsfonts,amsmath,amssymb}
\usepackage{fixltx2e} 
\usepackage{booktabs}

\title{\bf A Modified EMD Algorithm and its Applications}           

\author{
{\bfseries Mayer Humi$^1$}\\
$^1$Department of Mathematical Sciences, \\
Worcester Polytechnic Institute, \\
Worcester, MA 01609, USA 
}

\begin{document}

\maketitle                        

\begin{abstract}

The classical EMD algorithm has been used extensively in the literature to 
decompose signals that contain nonlinear waves. However when a signal contain 
two or more frequencies that are close to one another the decomposition might 
fail. In this paper we propose a new formulation of this algorithm which is 
based on the zero crossings of the signal and show that it performs well 
even when the classical algorithm fail. We address also
the filtering properties and convergence rate of the new algorithm versus 
the classical EMD algorithm. These properties are compared then 
to those of the principal component algorithm (PCA). Finally we apply this 
algorithm to the detection of gravity waves in the atmosphere. 

\end{abstract}

\vspace{1em}

\section{Introduction}

In scientific literature there exist many classical sets of functions 
which can decompose a signal in terms of "simple" functions. For 
example Taylor or Fourier expansions are used routinely in scientific and 
engineering applications.(and many other exist). However in all these
expansions the underlying functions are not intrinsic to the signal itself
and a precise approximation to the original signal might require a large
number of terms. This problem become even more acute when the signal is 
non-stationary and the process it represents is nonlinear.

To overcome this problem many researchers used in the past 
the "principal component algorithm" (PCA) to come up with an
"adaptive" set of functions which approximate a given signal. 
A new approach to this problem emerged in the late 1990's when a NASA team 
has developed the "Empirical Mode Decomposition" algorithm(EMD) which attempt 
to decompose a signal in terms of it "intrinsic mode functions"(IMF) through
a "sifting algorithm". A patent for this algorithm has been issued [1]. 

The EMD algorithm is based on the following quote [2]:
"According to Drazin the first step of data analysis is to examine
the data by eye. From this examination, one can immediately identify
the different scales directly in two ways: by the time lapse between
successive alterations of local maxima and minima and by the time lapse
between the successive zero crossings....We have decided to to adopt the
time lapse between successive extrema as the definition of the time
scale for the intrinsic oscillatory mode" 

A step by step description of the EMD sifting algorithm is as follows:
\begin{enumerate}
\item Let be given a function $f(t)$ which is sampled at discrete times 
$\{t_k, k=1,\ldots n\}$.
\item let $h_0(k)=f(t_k)$.
\item Identify the max and min of $h_0(k)$.
\item Create the cubic spline curve $M_x$ that connects the maxima points.
Do the same for the minima $M_n$. This creates an envelope for
$h_0(k)$.
\item At each time $t_k$ evaluate the mean $m_k$ of $M_x$ and $M_n$
($m_k$ is referred to as the sifting function).
\item Evaluate $h_1(k)=h_0(k)-m_k$.
\item If norm of $||h_0-h_1|| < \epsilon$ for some predetermined 
$\epsilon$ set the first intrinsic function $IMF_1=h_1$ (and stop).
\item if the criteria of (7) are not satisfied  set $h_0(k)=h_1(k)$
and return to (3) ("Sifting process"). 
\end{enumerate}

The algorithm has been applied successfully in various physical 
applications. However as has been observed by Flandrin [3] and others
the EMD algorithm fails in many cases where the data contains two
or more frequencies which are close to each other.

To overcome this difficulty we propose hereby a modification of the EMD
algorithm by replacing steps $4$ and $5$ in the description above by 
the following:

4. find the midpoints between two consecutive maxima and minima
and let $N_k$ be the values of $h_0$ at these points. \\
5. Create the spline curve $m_k$ that connects the points $N_k$.

The essence of this modification is the replacement of the mean 
which is evaluated by the EMD algorithm as the average of the max-min
envelopes by the spline curve of the mid-points between the maxima and minima.
This is in line with the observation by Drazin (which was referred to above)
that the scales inherent to the data can be educed either from the max-min 
or its zero crossing. In the algorithm we propose hereby we mimic the 
"zero-crossings" by the mid-points between the max-min. 
 
It is our objective in this paper to justify this modification of the EMD 
algorithm through some examples and theoretical work. The plan of the
paper is as follows: In Sec. $2$ we provides examples of signals composed
two or three close frequencies (with and without noise) where the classical 
EMD algorithm fails but the modified one yields satisfactory results.
In Sec. $3$ we carry out analytical analysis of the two algorithms
which are applied to the same signal. In Sec. $4$ we discuss the convergence
rate, resolution and related issues concerning the classical and new 
"midpoint algorithm" . Sec. $5$ address the application of this algorithm 
to atmospheric data and in Sec. $6$ we compare the EMD and PCA algorithms

\setcounter{equation}{0}
\section{Examples and Comparisons}

Extensive experimentations were made to test and verify the efficiency of 
the modified algorithm. We present here the results of one of these tests 
in which the signal contains three close frequencies. (In our tests we 
considered also the effects of noise and phase shifts among the different 
frequencies)

\begin{equation}
\label{2.1}
f(t)= \frac{1}{3}[\cos(\omega_1 t)+\cos(\omega_2 t)+\cos(\omega_3 t)]
\end{equation} 
where
$$
\omega_1=12\omega_0,\,\,\, \omega_2=10\omega_0,\,\,\, \omega_3=8\omega_0,
\,\,\, \omega_0=\frac{\pi}{256}.
$$
To apply the EMD algorithm to this signal, we used a discrete representation 
of it over the interval $[-2048,2048]$ by letting $t_{k+1}-t_{k}=1$,
$k=1,\ldots,4097$.

The results of the signal decompositions into IMFs and a comparison
these IMFs with the frequencies present in the original signal are presented 
in figures $1-5$. In all these figures the red lines represent the 
frequencies in the original signal (or its power spectrum) 
and the blue lines the corresponding intrinsic mode functions or their
power spectrum which were obtained by the midpoint algorithm. 

Fig. $1$ is a plot of the data for the signal described by (\ref{2.1}).
Fig. $2$ represents the first IMF in the decomposition (versus 
the leading frequency in the data) while Figs. $3-5$ depict the spectral 
density distribution for the first three IMFs versus those related to 
the original frequencies in the data. It should be observed that although 
the amplitude of the spectral densities in these plots are different 
(especially for IMF 3) the maxima of the spectral density in each plot is 
very close to the original one.

The EMD algorithm is a high pass filter. For the $n-th$ iteration of the 
filter its efficiency is measured by the parameter $\alpha$ which is 
defined by
$$
Y_n= \alpha_n Y_{n-1} + \alpha(X_n-X_{n-1})
$$
where $X_k$ and $Y_{k}$ are the input and output of the $k-th$ iteration. 
Fig $6$ present the value of the parameter $\alpha$ as a function of the
iteration number for first IMF derived from the data of the signal
in (\ref{2.1}).   

\setcounter{equation}{0}
\section{Some Analytical Insights}

To obtain analytical insights about the performance of the 
EMD-midpoint algorithm we considered the following signal
\begin{equation}
\label{3.1}
f(t)=\frac{1}{2}[\cos(\omega_4 t)+\cos(\omega_5t)],\,\,\,
\omega_4=\frac{3\pi}{64},\,\,\,\omega_5=\frac{\pi}{32}.
\end{equation} 
Since the ratio of the frequencies in this signal is a rational number 
the signal is actually periodic with period $p=128$ 
(See Fig. $7$) and the behavior of 
the classical versus the mid-point algorithm can be delineated 
analytically (i.e without discretizations). 

On the interval $[0,p]$ the extrema of the signal are given by
$\frac{df}{dt}=0$ and therefore it is easy to construct the 
spline approximation $S_{max}(t)$, $S_{min}(t)$ to the maximum and 
minimum points and compute their average. Similarly we can 
find the midpoints between the maxima and minima and evaluate the 
corresponding spline approximation $S_{mid}(t)$ to the signal at these points.
after one iteration of the sifting process the "sifted signal"
is given respectively by
\begin{equation}
\label{3.2}
h_{mn}(t)=f(t)-\frac{S_{max}(t)+S_{min}(t)}{2},
\end{equation} 
and
\begin{equation}
\label{3.3}
h_{mid}(t)=f(t)-S_{mid}(t).
\end{equation} 
The efficiency of the two algorithm can be deduced by projecting
these new signals on the Fourier components of the original signal. To
this end we compute
\begin{equation}
\label{3.4}
a_{mn}=\displaystyle\int_{0}^{p}h_{mn}(t)\cos(\omega_4 t) dt,\,\,\
b_{mn}=\displaystyle\int_{0}^{p}h_{mn}(t)\sin(\omega_4 t)dt.
\end{equation} 
\begin{equation}
\label{3.5}
c_{mn}=\displaystyle\int_{0}^{p}h_{mn}(t)\cos(\omega_5 t) dt,\,\,\
d_{mn}=\displaystyle\int_{0}^{p}h_{mn}(t)\sin(\omega_5 t)dt.
\end{equation} 
and
\begin{equation}
\label{3.6}
a_{mid}=\displaystyle\int_{0}^{p}h_{mid}(t)\cos(\omega_4 t) dt,\,\,\
b_{mid}=\displaystyle\int_{0}^{p}h_{mid}(t)\sin(\omega_5 t)dt.
\end{equation} 
\begin{equation}
\label{3.7}
c_{mid}=\displaystyle\int_{0}^{p}h_{mid}(t)\cos(\omega_4 t) dt,\,\,\
d_{mid}=\displaystyle\int_{0}^{p}h_{mid}(t)\sin(\omega_5 t)dt.
\end{equation} 
The amplitude of the Fourier components of the two frequencies in the 
classical EMD algorithm is 
\begin{equation}
\label{3.8}
A_{mn}=\sqrt{a_{mn}^2+b_{mn}^2},\,\,\,\, B_{mn}=\sqrt{c_{mn}^2+d_{mn}^2}.
\end{equation} 
Similarly for the mid-point algorithm we
\begin{equation}
\label{3.9}
A_{mid}=\sqrt{a_{mid}^2+b_{mid}^2},\,\,\,\, B_{mid}=\sqrt{c_{mid}^2+d_{mid}^2}.
\end{equation} 
The objective of the sifting process is to eliminate one of the
Fourier components in favor of the other. As a result the first IMF 
will contains, upon convergence, only one of the Fourier components in the 
original signal. Therefore the efficiency of the two algorithm can be 
inferred by comparing $A_{mn}$ versus $B_{mn}$ and $A_{mid}$ versus
$B_{mid}$. Computing the integrals that appear in eqs.(\ref{3.4})-(\ref{3.7})
we obtain
\begin{equation}
\label{3.11}
A_{mn}=31.63346911,\,\,\,\, B_{mn}=29.70292046,
\end{equation}
\begin{equation}
\label{3.12}
A_{mid}=34.19647843,\,\,\,\, B_{mid}=20.81145369.
\end{equation} 
These results show that after one iteration the classical EMD did not 
separate the two frequencies effectively. On the other hand the mid-point 
algorithm performed well.

\setcounter{equation}{0}
\section{Convergence Rates}

To compare the convergence rates of the classical versus the midpoint 
algorithm we considered three cases all of which were composed of two 
frequencies. In the first case the two frequencies were well separated.
In the second case the two frequencies were close while in the third case
they were almost "overlapping". In all cases the signal was given by
$$
f(t)=\frac{1}{2}(\cos\omega_1 t +cos\omega_2 t)
$$
This signal was discretized on the interval $[-2048,2048]$ with $\Delta t=1$.
 
For the first case the two frequencies were
$$
\omega_1=12\omega,\,\,\,\, \omega_2=8\omega,\,\,\,\, \omega=\frac{\pi}{256}.
$$
As can be expected both the classical and midpoint algorithm were able to
discern the individual frequencies through the sifting algorithm. However 
it took the classical algorithm $59$ iterations to converge to the first IMF. 
On the other hand the midpoint algorithm converged in only $7$ iterations 
(using the same convergence criteria). We wish to point out also that the
midpoint algorithm has a lower computational cost than the classical 
algorithm. It requires in each iteration the computation of only one spline 
interpolating polynomial. On the other hand the classical algorithm requires 
two such polynomials, one for the maximum points and one for the minimum points.

For the second test the frequencies were
$$
\omega_1=\frac{\pi}{24}+\frac{\pi}{288},\,\,\,\,
\omega_2=\frac{\pi}{24}-\frac{\pi}{288}
$$
that is the difference between the two frequencies is $\frac{\pi}{144}$.

In this case the midpoint algorithm was able to separate the two frequencies.
Fig $8$ and Fig $9$ compare the power spectrum of the original frequencies 
versus those of $IMF_1$ and $IMF_2$ which were obtained through this algorithm. 
Convergence to $IMF_1$ was  obtained in 18 iterations and $IMF_2$ was obtained 
by $7$ additional iterations.

The classical EMD algorithm did converge to $IMF_1$ in $45$ iterations but 
the power spectrum of this $IMF$ deviated significantly from the first 
frequency in the signal(See Fig $10$). $IMF_2$ failed (completely) to 
detected correctly the second frequency.

In third case the frequencies were
$$
\omega_1=\frac{\pi}{24}+\frac{\pi}{1000},\,\,\,\,
\omega_1=\frac{\pi}{24}-\frac{\pi}{1000}.
$$
In this case the classical algorithm was unable to separate the two 
frequencies i.e $IMF_1$ contained both frequencies (See Fig $11$). 
The midpoint algorithm did somewhat better but the resolution was not 
complete (See Fig $12$). Moreover the sifting process in both cases led 
to the creation of "ghost frequencies" which were not present in the 
original signal. 

At this juncture one might wonder if a "hybrid algorithm" whereby the sifting
function is the average (or some similar combination) of those obtained by 
the classical and midpoint algorithms might outperform the separate algorithms
(in spite of the obvious additional computational cost). However our 
experimentations with such algorithm did not yield the desired results
(i.e. the convergence rate and resolution did not improve).

\setcounter{equation}{0}
\section{Applications to Atmospheric Data}

There have been recent interest in the observation and properties of gravity
waves which are generated when wind is blowing over terrain. In part this 
interest stems from the fact that these waves carry energy and accurate 
measure of this data is needed to improve the performance of numerical 
weather prediction models.

As part of this scientific campaign the USAF flew several balloons that 
collected information about the pressure and temperature as a function of
height. The temperature data collected by one of these balloons is presented 
in Fig. $13$ [6]. To analyze this signal we detrended first it by subtracting 
its mean from the data. When the mid-point EMD algorithm was applied to this 
detrended-signal the first IMF extracted the experimental noise from 
while the second and third IMFs educed clearly the gravity waves  
(the second IMF is depicted in Fig. $14$). On the other hand the classical
EMD algorithm failed to educe these waves from the detrended-signal.

Subtracting the gravity waves that were detected by the mid-point algorithm 
from the detrended-signal we obtain the "turbulent residuals" whose spectrum
is shown in Fig $15$. The slope of this signal in the "inertial frequency 
range" is $-2.7$ which corresponds well with the fact that the flow in
stratosphere is "quasi two-dimensional" [7-9].

\setcounter{equation}{0}
\section{EMD or PCA- A Comparison}

Before the emergence of the EMD algorithm an adaptive data analysis 
was provided by the "Principal Component Algorithm"(PCA) which is referred
to also as the "Karahunan-Loeve (K-L) decomposition algorithm".
(For a review see [10]) Here we shall give only a brief overview
of this algorithm within in the geophysical context.

Let a signal be represented by a a time series  $X$  (of length  $N$) 
of some variable.We first determine a time delay  $\Delta$  for which the
points in the series are decorrelated.  Using  $\Delta$  we create
$n$  copies of the original series
$$
X(k),\;\;X(d + \Delta),\ldots,X(k + (n - 1)\Delta).
$$
(To create these one uses either periodicity or choose to consider
shorter time-series).  Then one computes the auto-covariance matrix
$R = (R_{ij})$
\begin{equation}
\label{Delta}
R_{ij} = \displaystyle\sum^N_{k=1} X(k + i\Delta)X(k+j\Delta).
\end{equation}
Let  $\lambda_0 > \lambda_1, \ldots, > \lambda_{n-1}$  be the
eigenvalues of  $R$  with their corresponding eigenvectors
$$
\phi^i = (\phi^i_0, \ldots, \phi^i_{n-1}),\;\;i = 0, \ldots, n - 1.
$$
The original time series  $X$  can be reconstructed then as
\begin{equation}
\label{phi}
X(j) = \displaystyle\sum^{n-1}_{k=0} a_k(j)\phi^k_0
\end{equation}
where
\begin{equation}
\label{sum}
a_k(j) = \displaystyle\frac{1}{n}\displaystyle\sum^{n-1}_{i=0}X(j +
i\Delta)\phi^k_i.
\end{equation}
The essence of the PCA  is based on the recognition that
if a large spectral gap exists after the first  $m_1$  eigenvalues of
$R$  then one can reconstruct the mean flow (or the large component
( of the data by using only the first  $m_1$  eigenfunctions in
(\ref{phi}).  A recent refinement of this procedure due to Ghil et al
([10]) is that the data corresponding to eigenvalues between  $m_1
+ 1$  and up to the point  $m_2$  where they start to form a
``continuum'' represent waves.  The location of  $m_2$ can be
ascertained further by applying the tests devised by Axford [11] and
Dewan [7].

Thus the original data can be decomposed into mean flow, waves and
residuals (i.e. data corresponding to eigenvalues  
$m_2 + 1,\ldots, n - 1$  which we wish to interpret at least partly 
as turbulent residuals).

The crucial step in this algorithm is the determination of the 
points $m_1$ and $m_2$ whose position has to ascertained by additional tests
whose results might be equivocal.

We applied this algorithm to the geophysical data described in Sec. $5.1$
with $\Delta=96$ and computed the resulting spectrum of the correlation 
matrix $R$. This spectrum is depicted in Fig. $16$ . Based on this spectrum
we choose $m_1=3$ and $m_2=11$ we obtain the corresponding wave 
component of the signal that is shown in Fig. $17$. 

We conclude that while the PCA algorithm provides an alternative to 
the EMD algorithm the determination of the cutoff points is murky in 
many cases. However it will be advantageous if one apply the two algorithms
in tandem in order to obtain a clear cut confirmation of the results.

\section*{References}

\begin{itemize}

\item[1] N. E. Huang - USA Patent $\# 6,311,130 B 1$ , Date Oct 30,2001

\item[2] N. E. Huang et all, ``The empirical mode decomposition and the Hilbert
spectrum for nonlinear and non-stationary time series analysis",
Proceedings of the Royal Society Vol. 454 pp.903-995 (1998)

\item[3] Gabriel Rilling and Patrick Flandrin,  ``One or Two Frequencies? The
Empirical Mode Decomposition Answers", IEEE Trans. Signal Analysis 
Vol. 56 pp.85-95 (2008).

\item[4] Zhaohua Wu and Norden E. Huang, ``On the Filtering Properties of 
the Empirical Mode Decomposition, Advances in Adaptive Data Analysis", 
Volume: 2, Issue: 4 pp. 397-414. (2010)

\item[5] Albert Ayenu-Prah and  Nii Attoh-Okine, ``A Criterion for Selecting 
Relevant Intrinsic Mode Functions in Empirical Mode Decomposition",
Advances in Adaptive Data Analysis, Vol. 2, Issue: 1(2010) pp. 1-24.

\item[6] George Jumper, ``Private communication" (2001)

\item[7] Dewan, E.M., ``On the nature of atmospheric waves and
turbulence, Radio Sci."  20, p. 1301-1307 (1985).

\item[8] Kraichnan, R., ``On Kolmogorov inertial-range theories", 
J. Fluid Mech., 62, p. 305-330 (1974).

\item[9] Lindborg, E., ``Can the atmospheric kinetic energy spectrum
be explained by two dimensional turbulence", J. Fluid Mech, 
388, p. 259-288 (1999).

\item[10] C. Penland, M. Ghil and K.M. Weickmann, ``Adaptive
filtering and maximum entropy spectra, with application to changes in
atmospheric angular momentum'', J. Geophys. Res., 96, 
22659-22671 (1991).

\item[11] D. N. Axford, ``Spectral analysis of aircraft observation
of gravity waves'', Q.J. Royal Met. Soc., 97, 313-321 (1971).
\end{itemize}

\includegraphics[scale=1,height=160mm,angle=0,width=180mm]{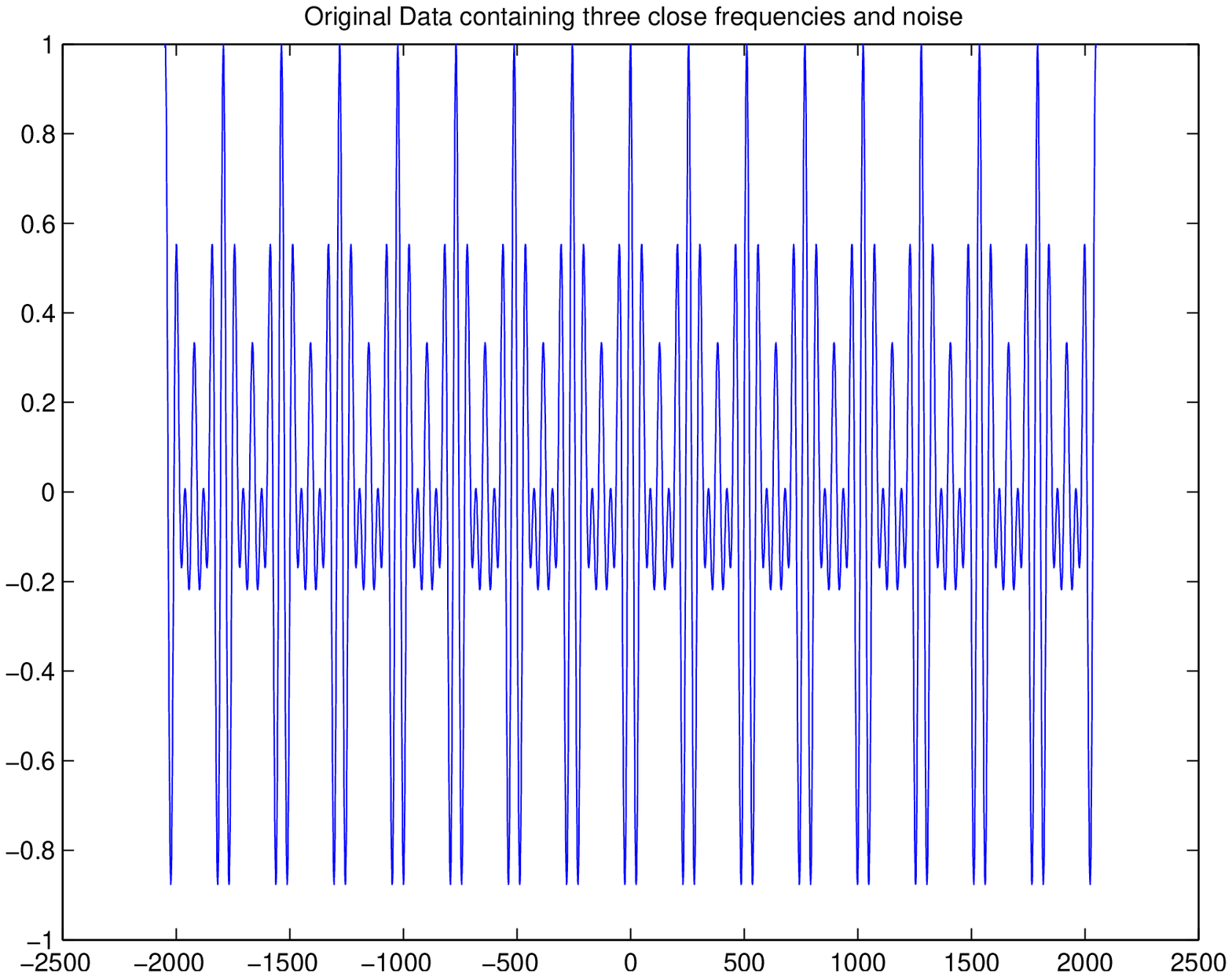}
\includegraphics[scale=1,height=160mm,angle=0,width=180mm]{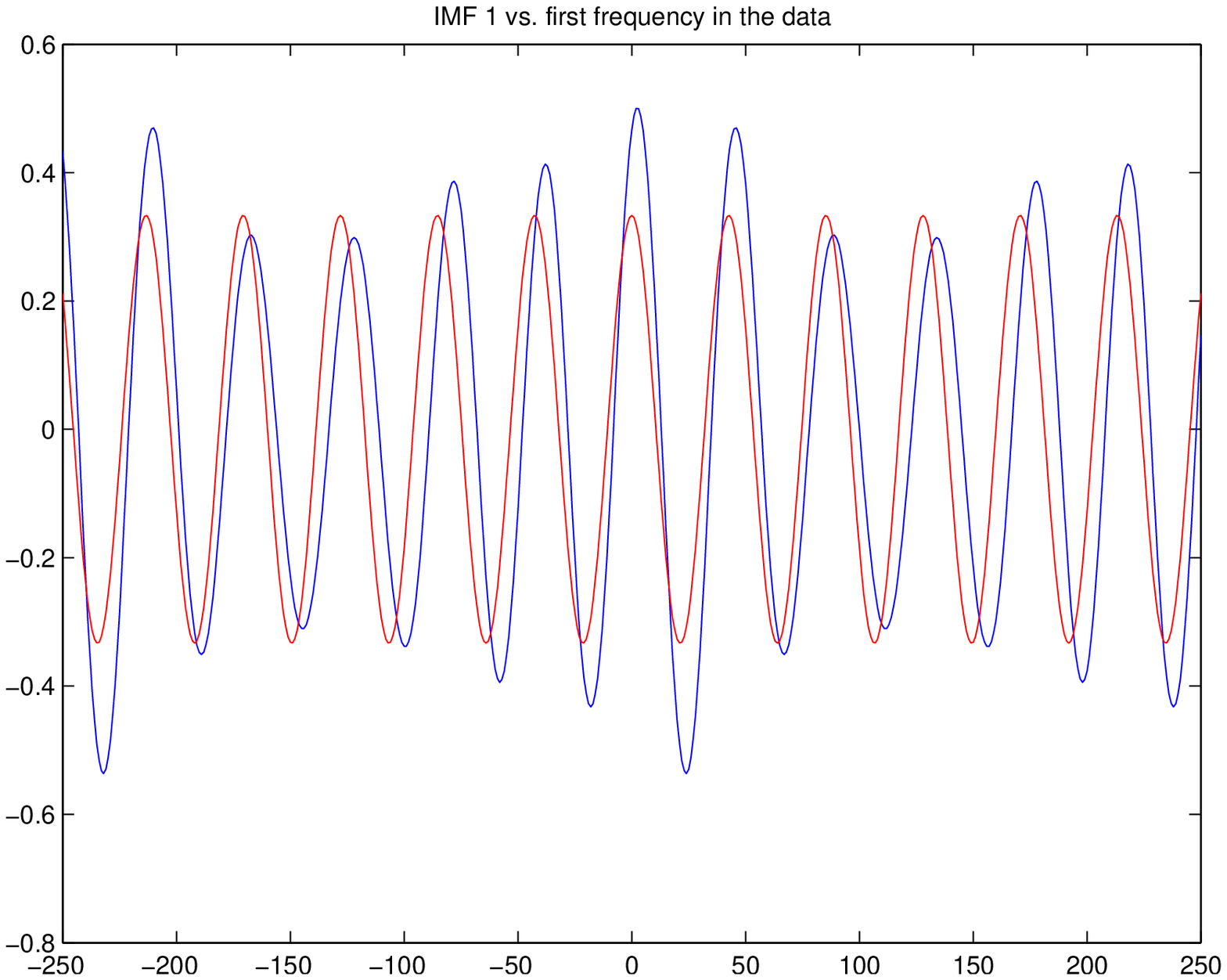}
\includegraphics[scale=1,height=160mm,angle=0,width=180mm]{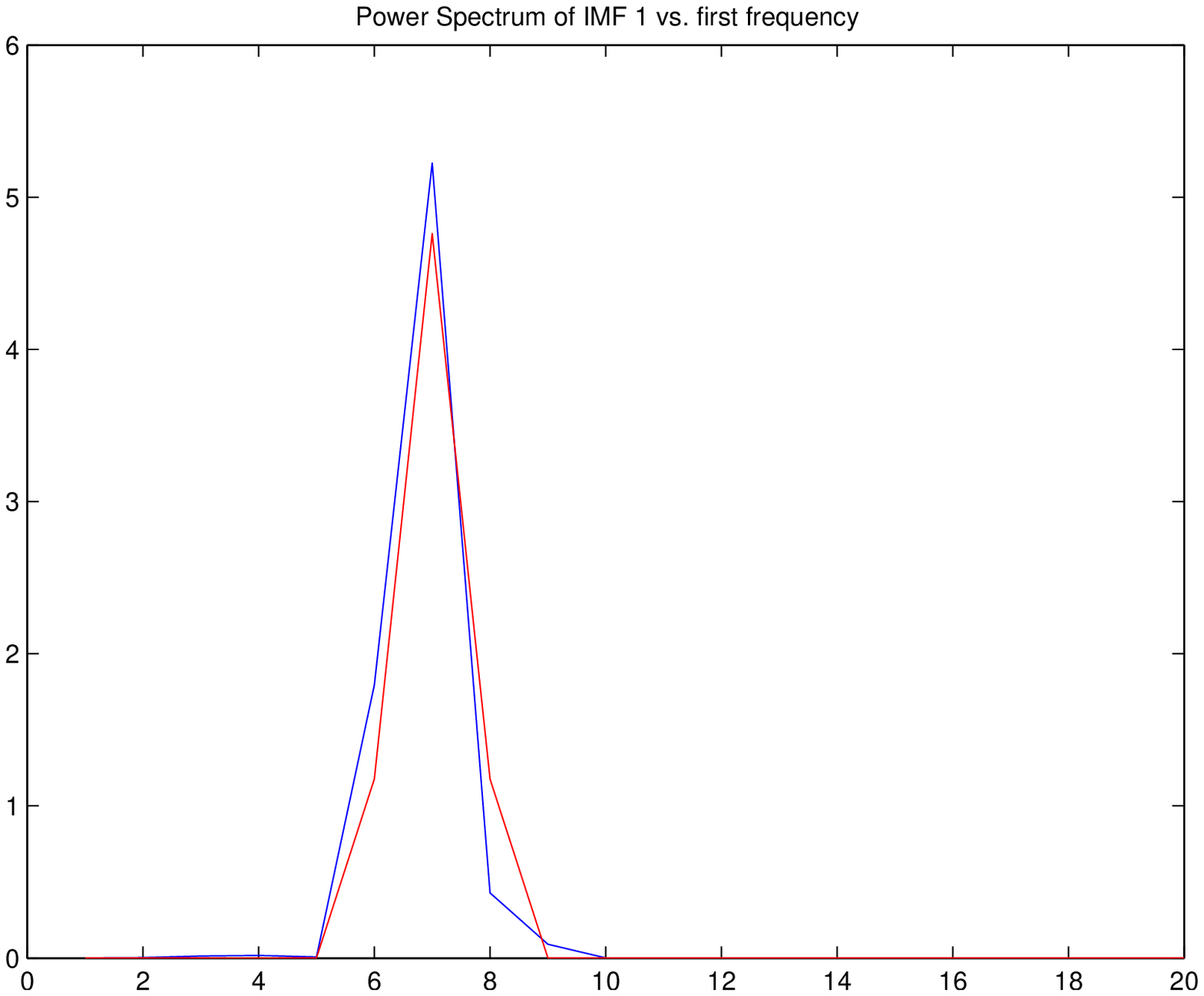}
\includegraphics[scale=1,height=160mm,angle=0,width=180mm]{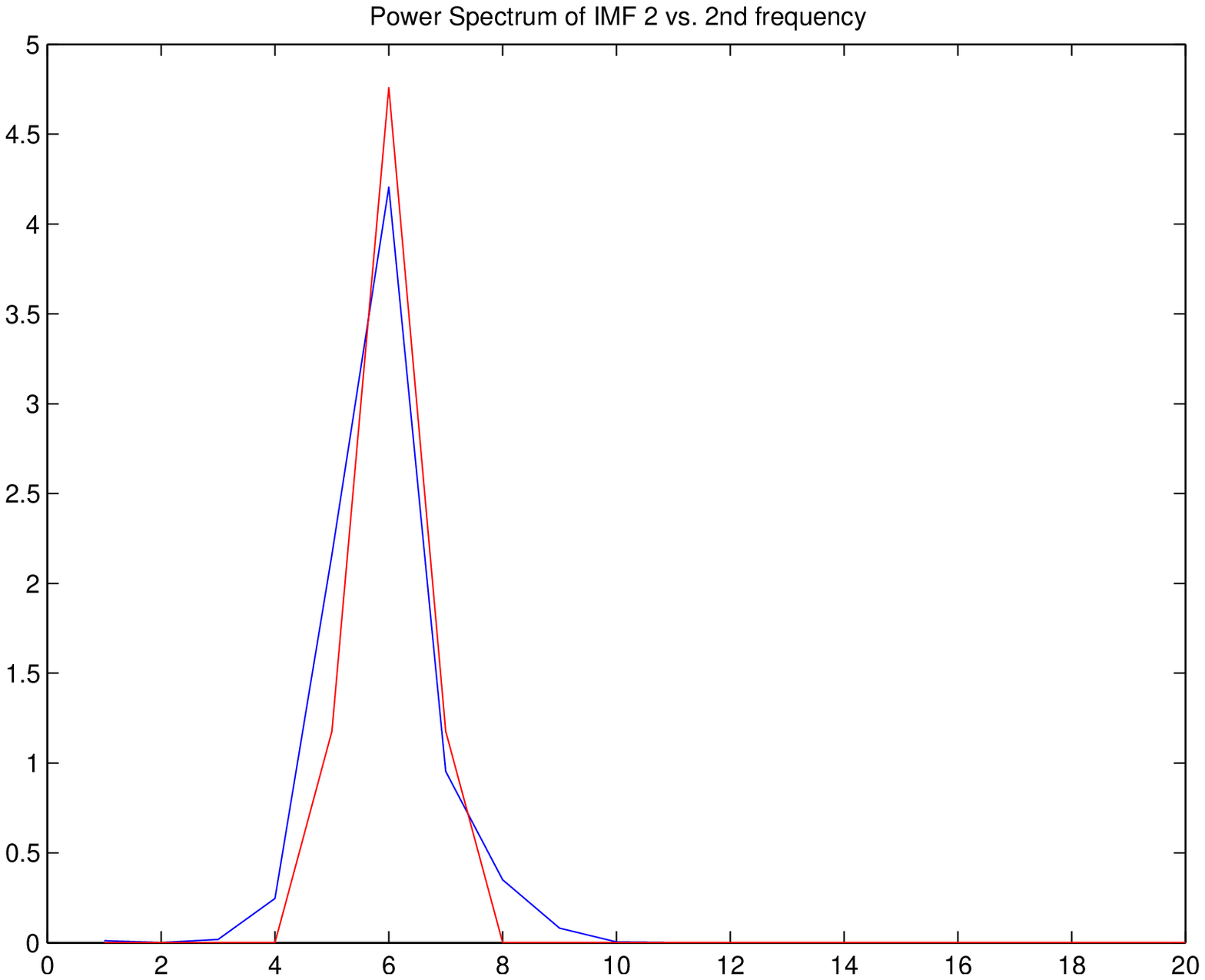}
\includegraphics[scale=1,height=160mm,angle=0,width=180mm]{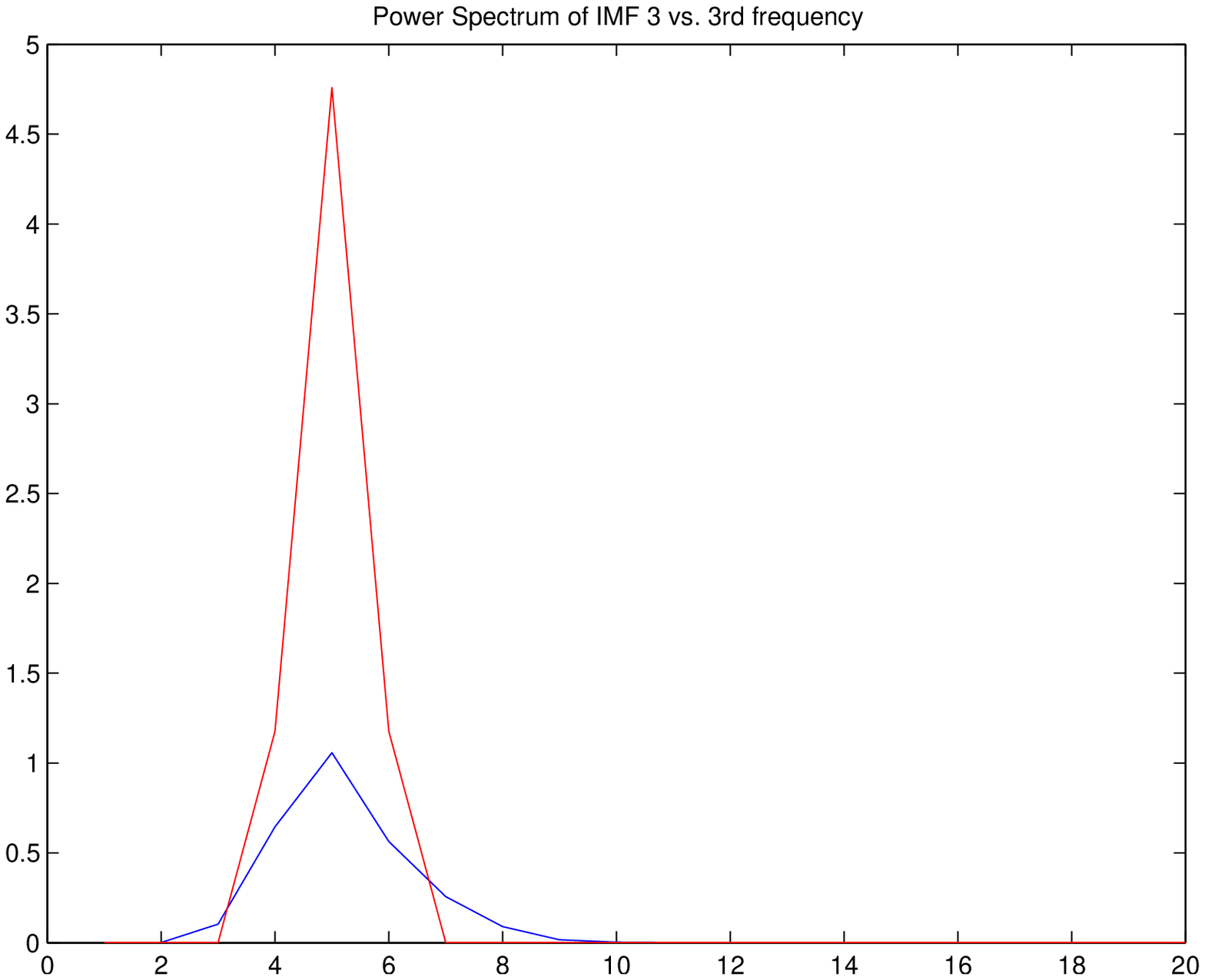}
\includegraphics[scale=1,height=160mm,angle=0,width=180mm]{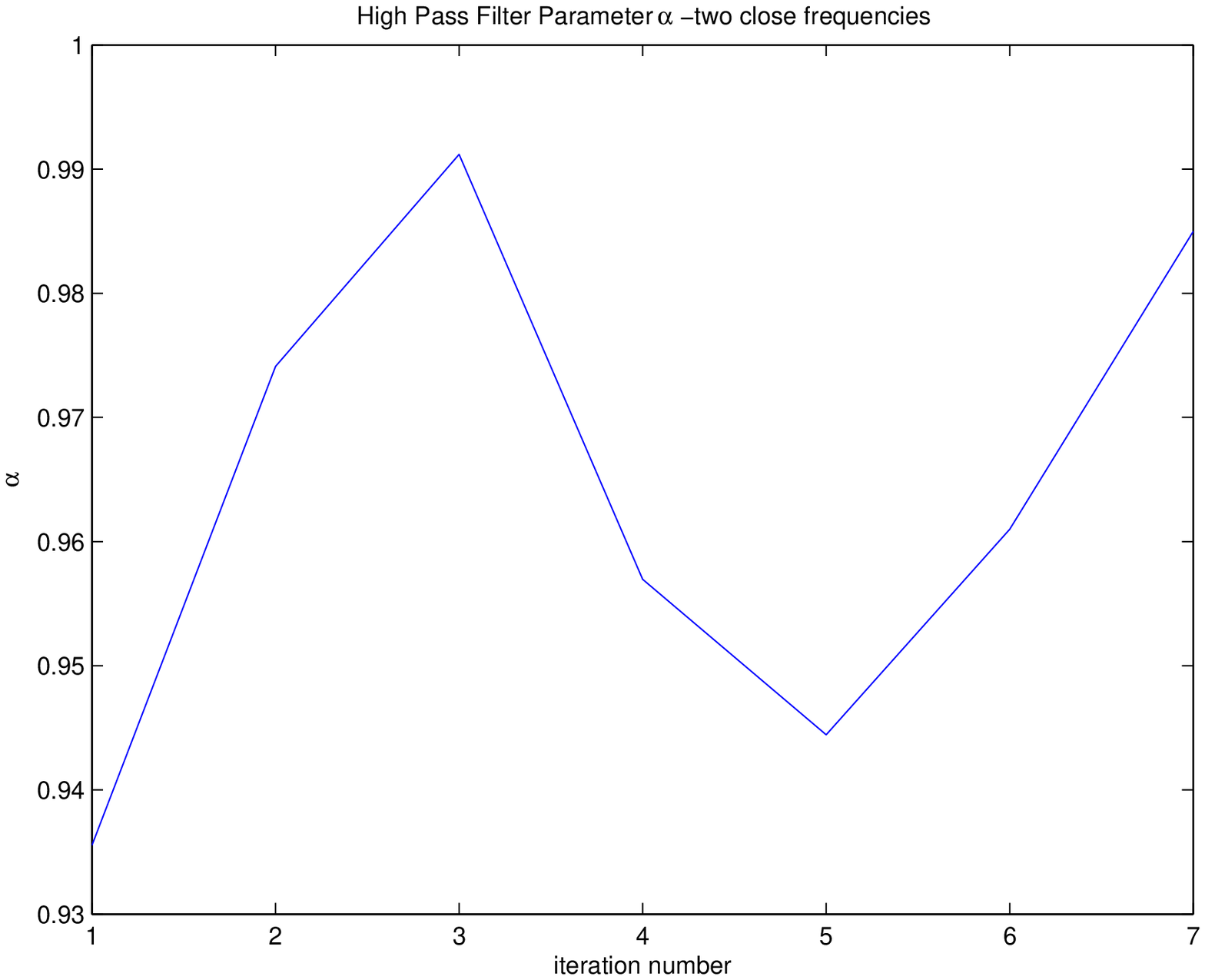}
\includegraphics[scale=1,height=160mm,angle=0,width=180mm]{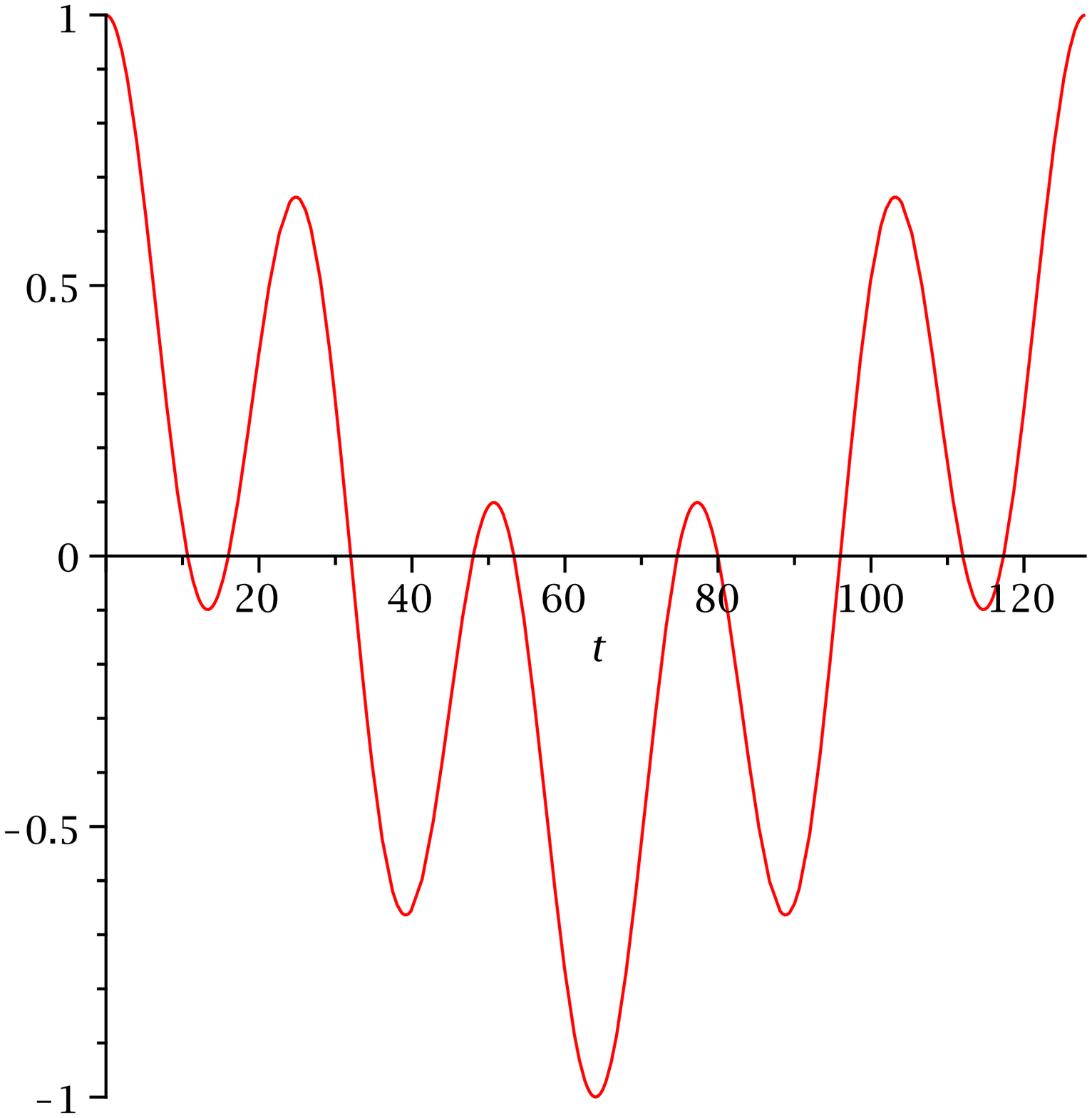}
\includegraphics[scale=1,height=160mm,angle=0,width=180mm]{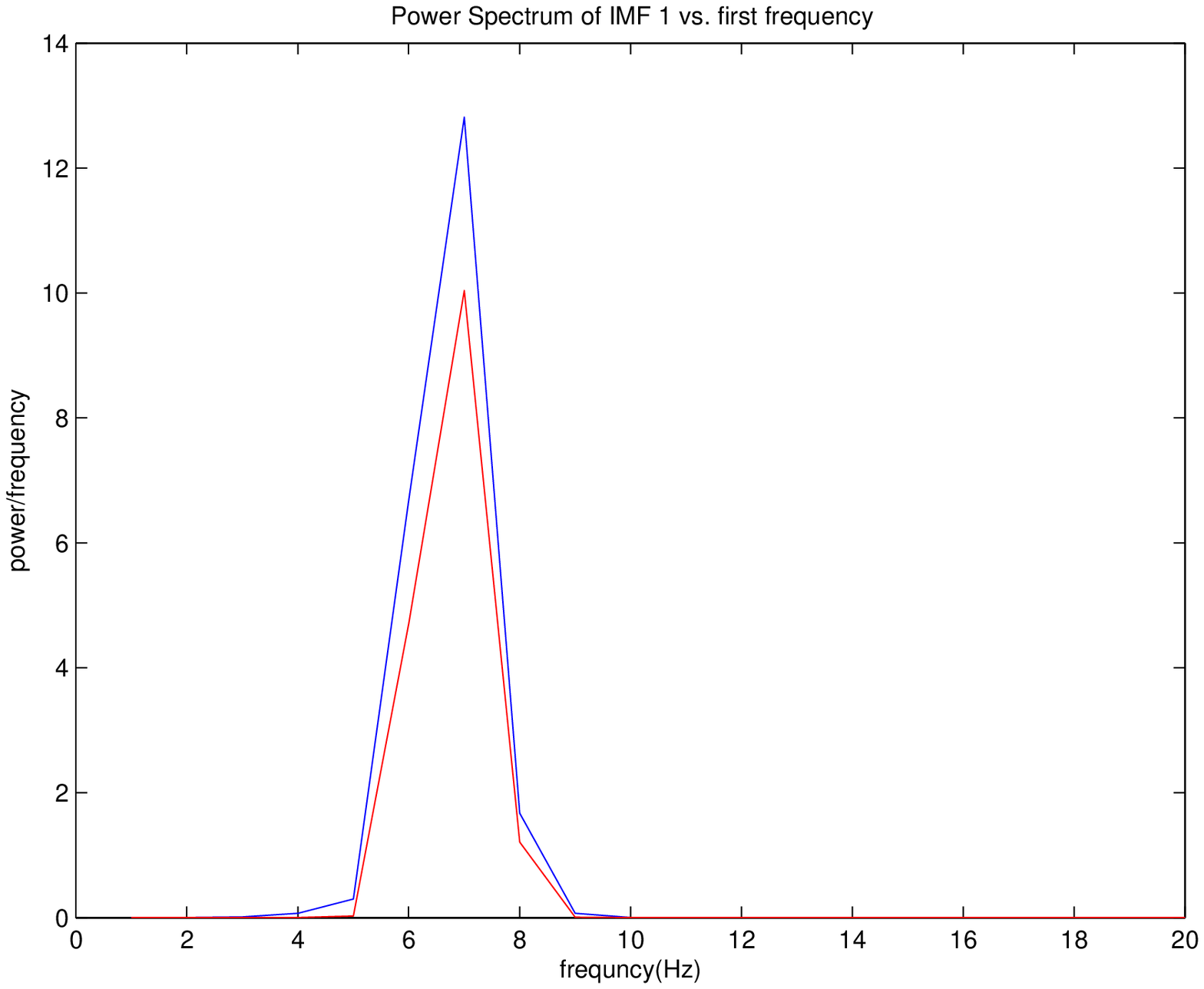}
\includegraphics[scale=1,height=160mm,angle=0,width=180mm]{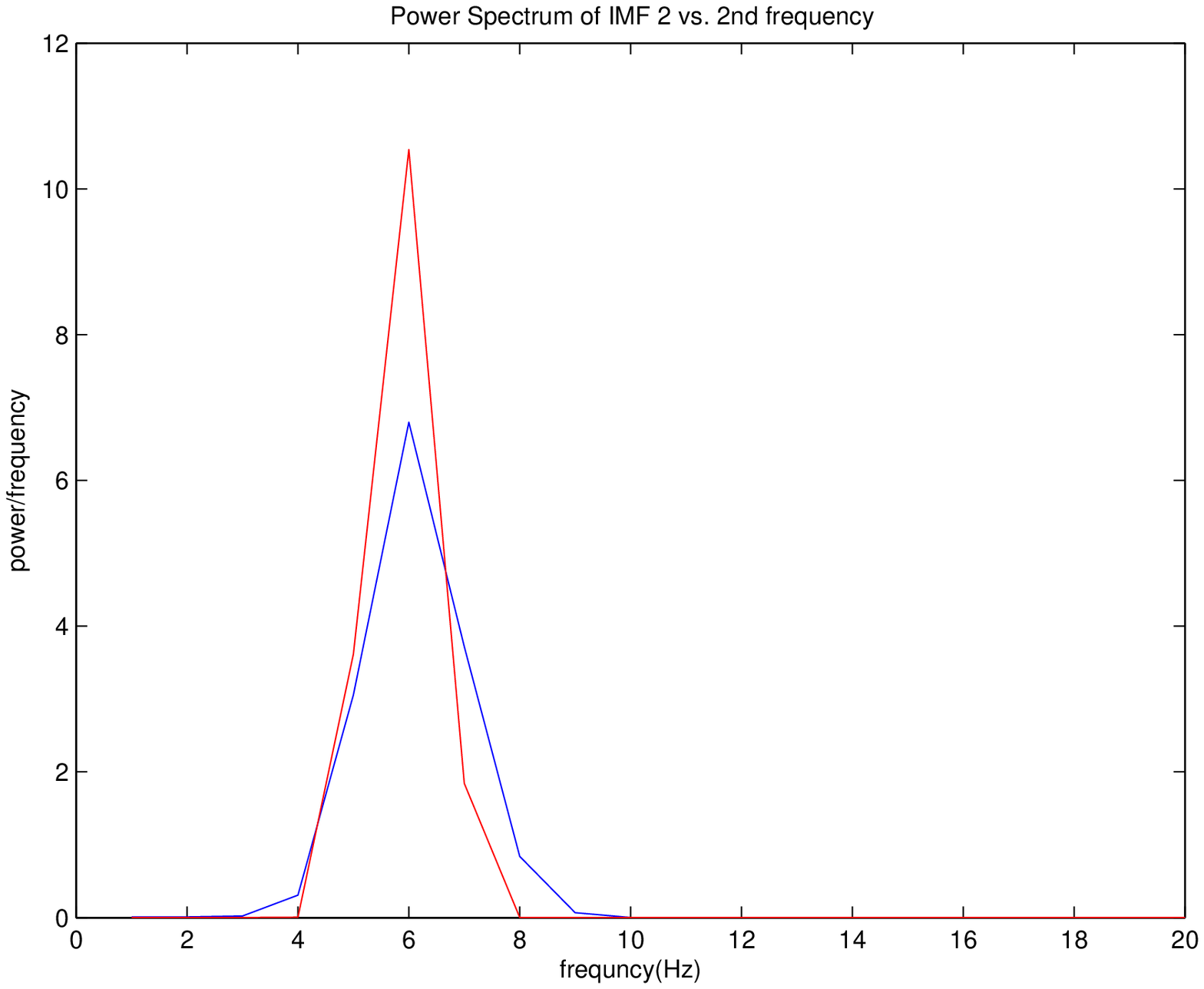}
\includegraphics[scale=1,height=160mm,angle=0,width=180mm]{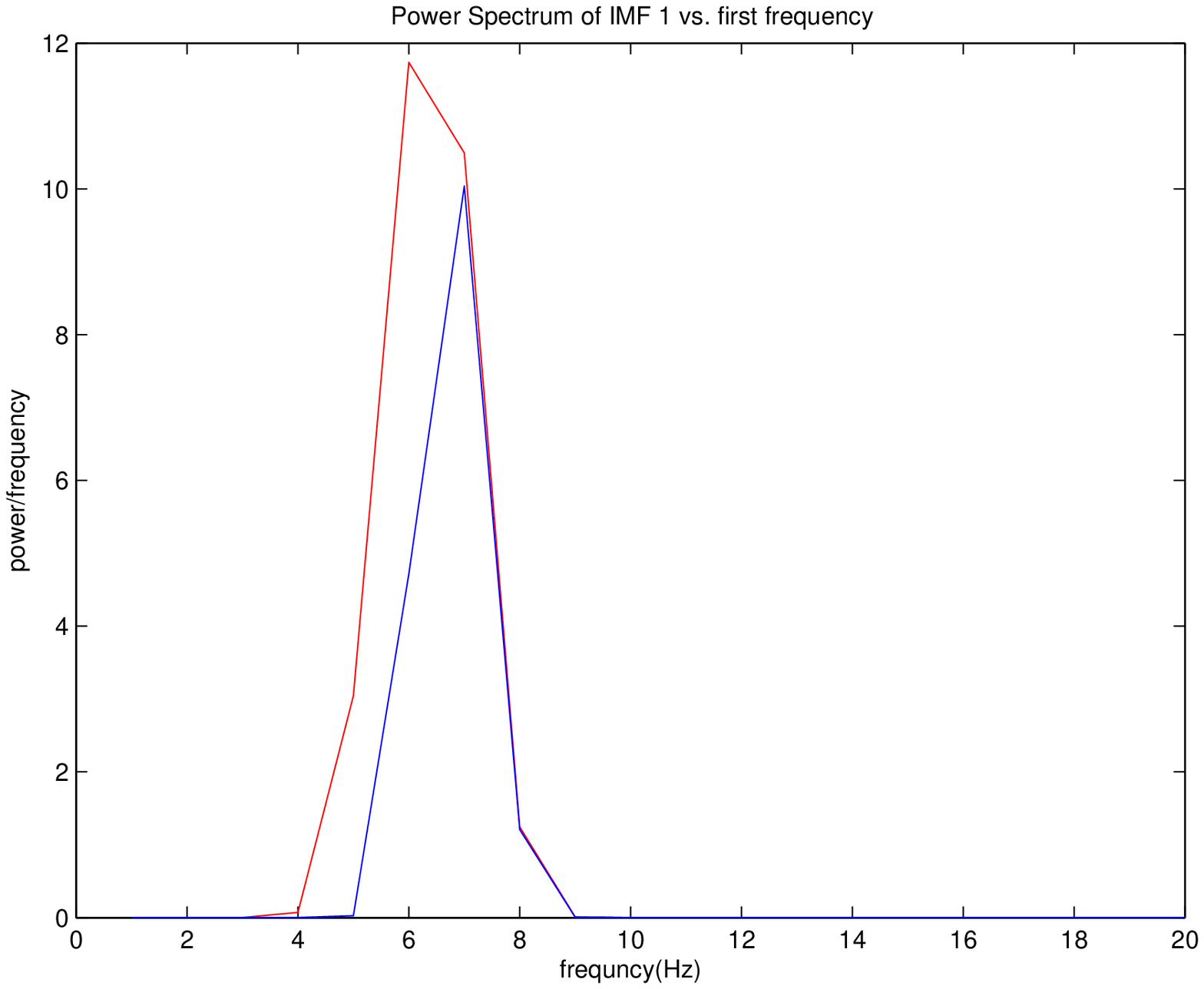}
\includegraphics[scale=1,height=160mm,angle=0,width=180mm]{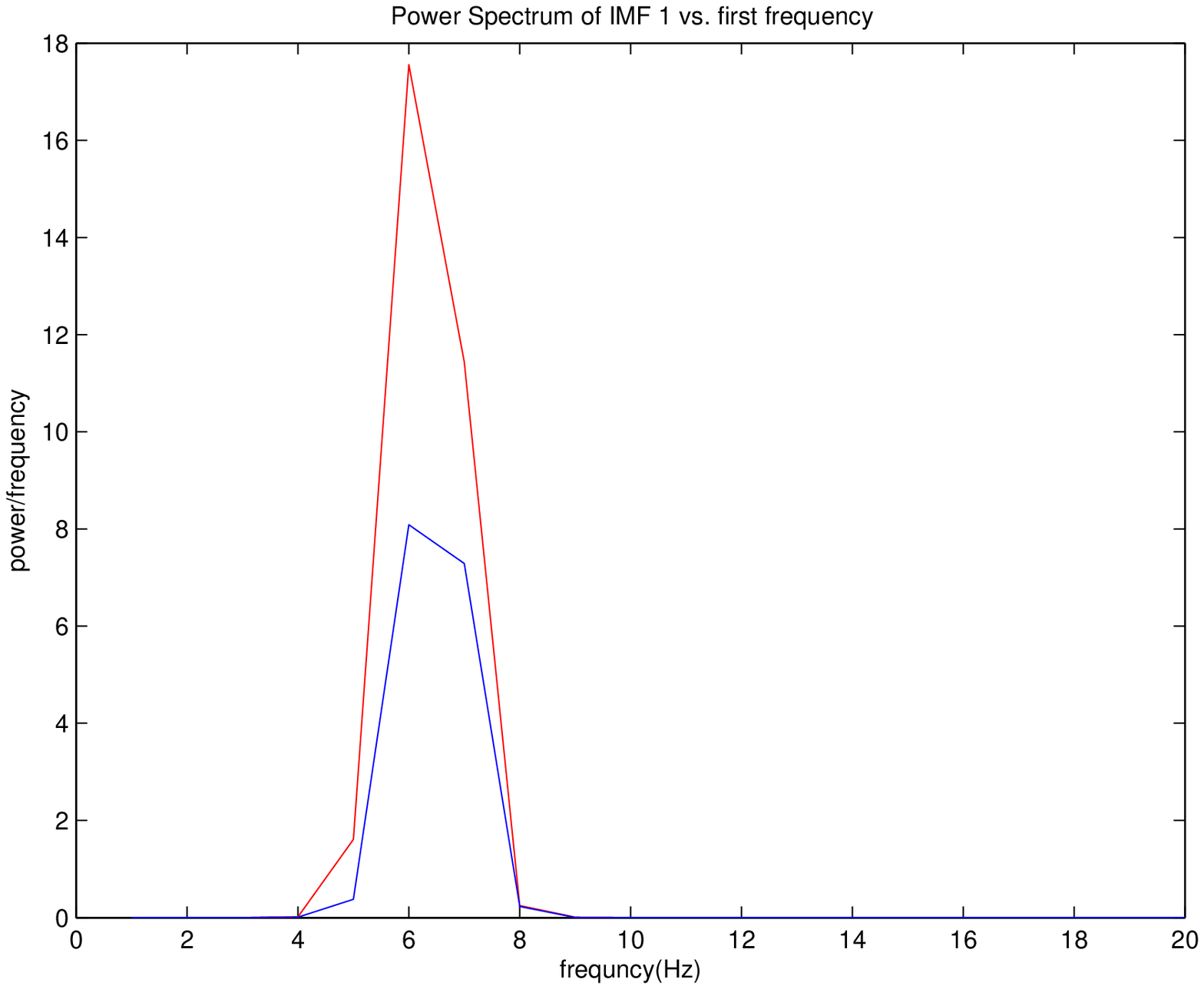}
\includegraphics[scale=1,height=160mm,angle=0,width=180mm]{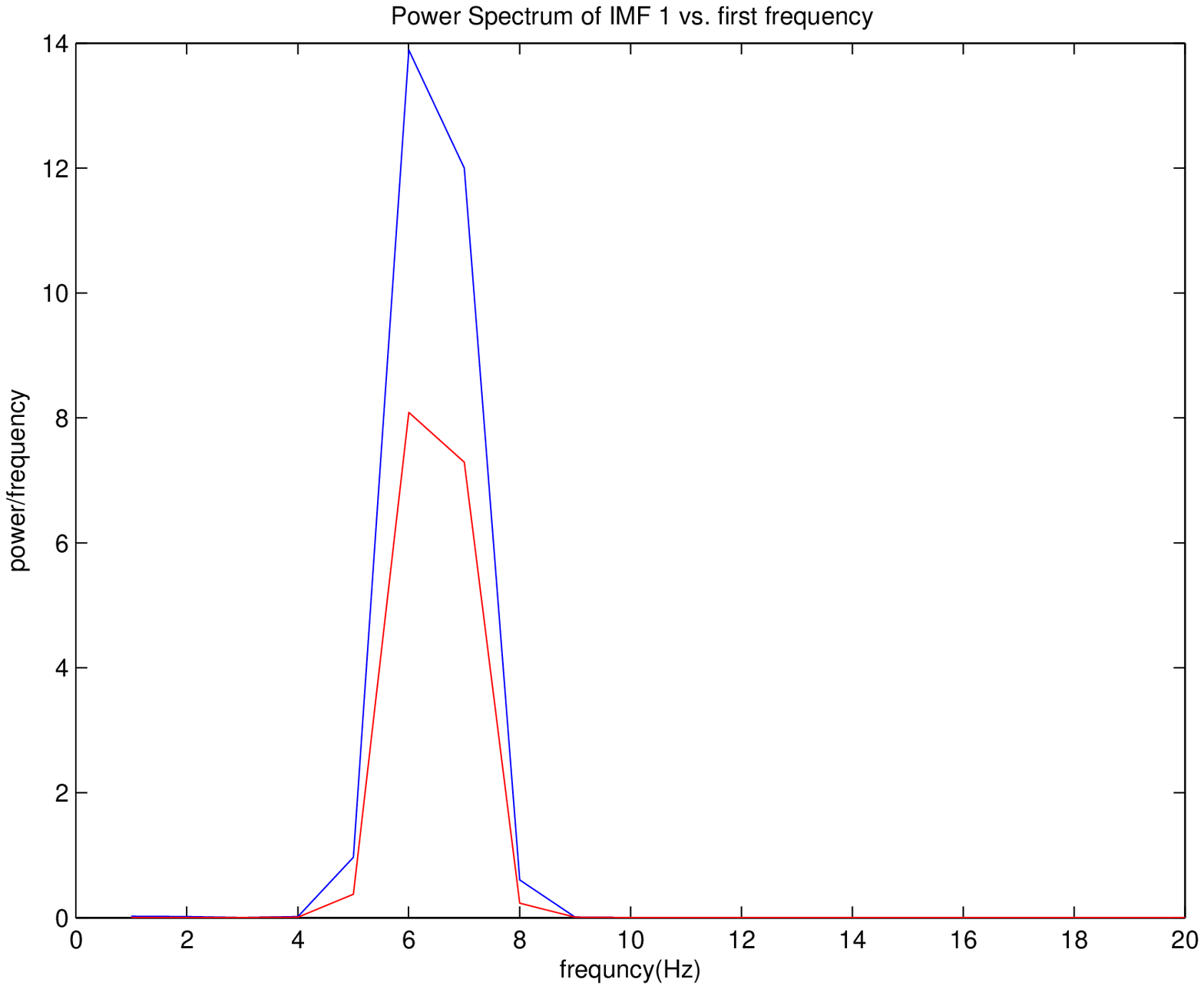}
\includegraphics[scale=1,height=160mm,angle=0,width=180mm]{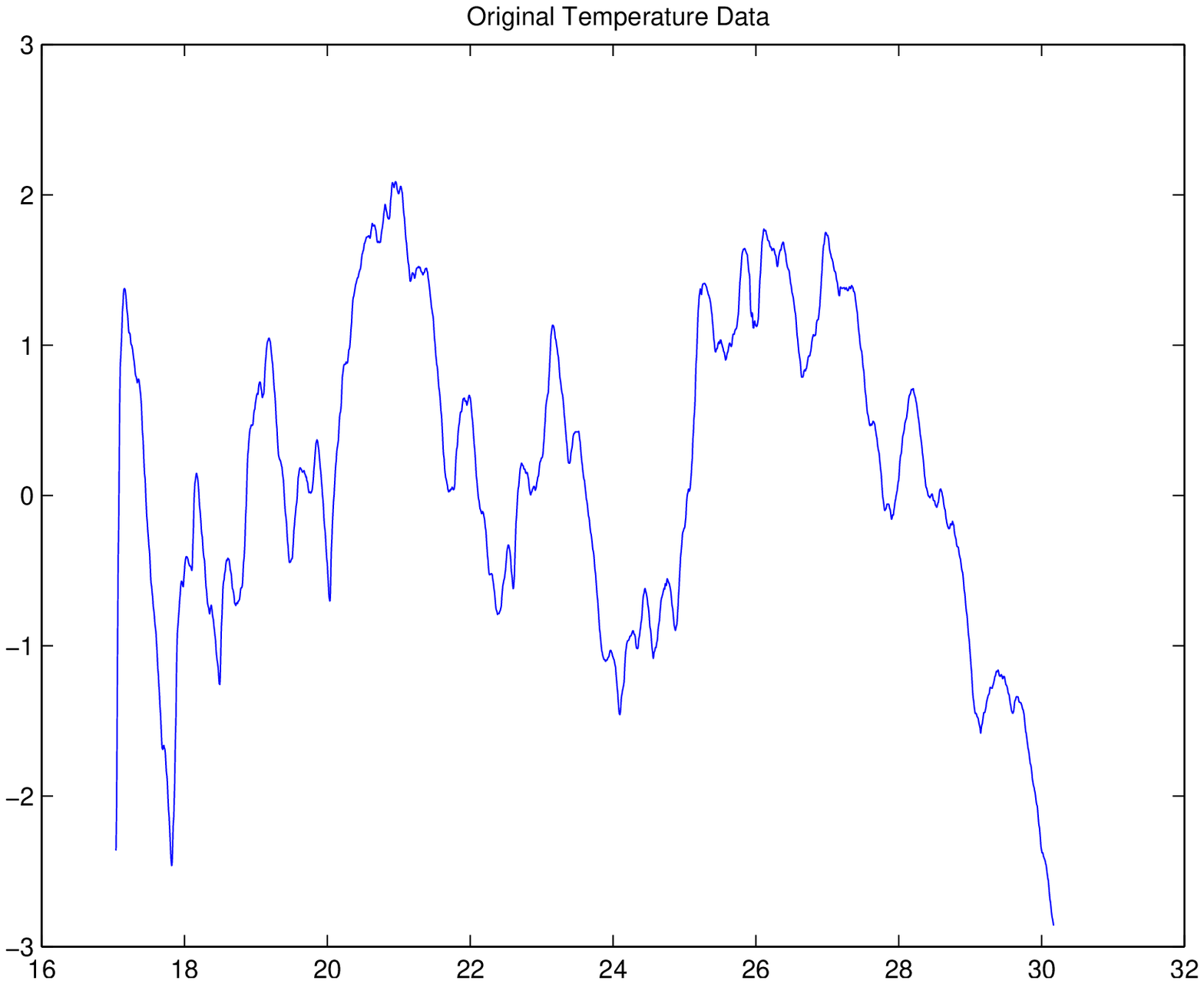}
\includegraphics[scale=1,height=160mm,angle=0,width=180mm]{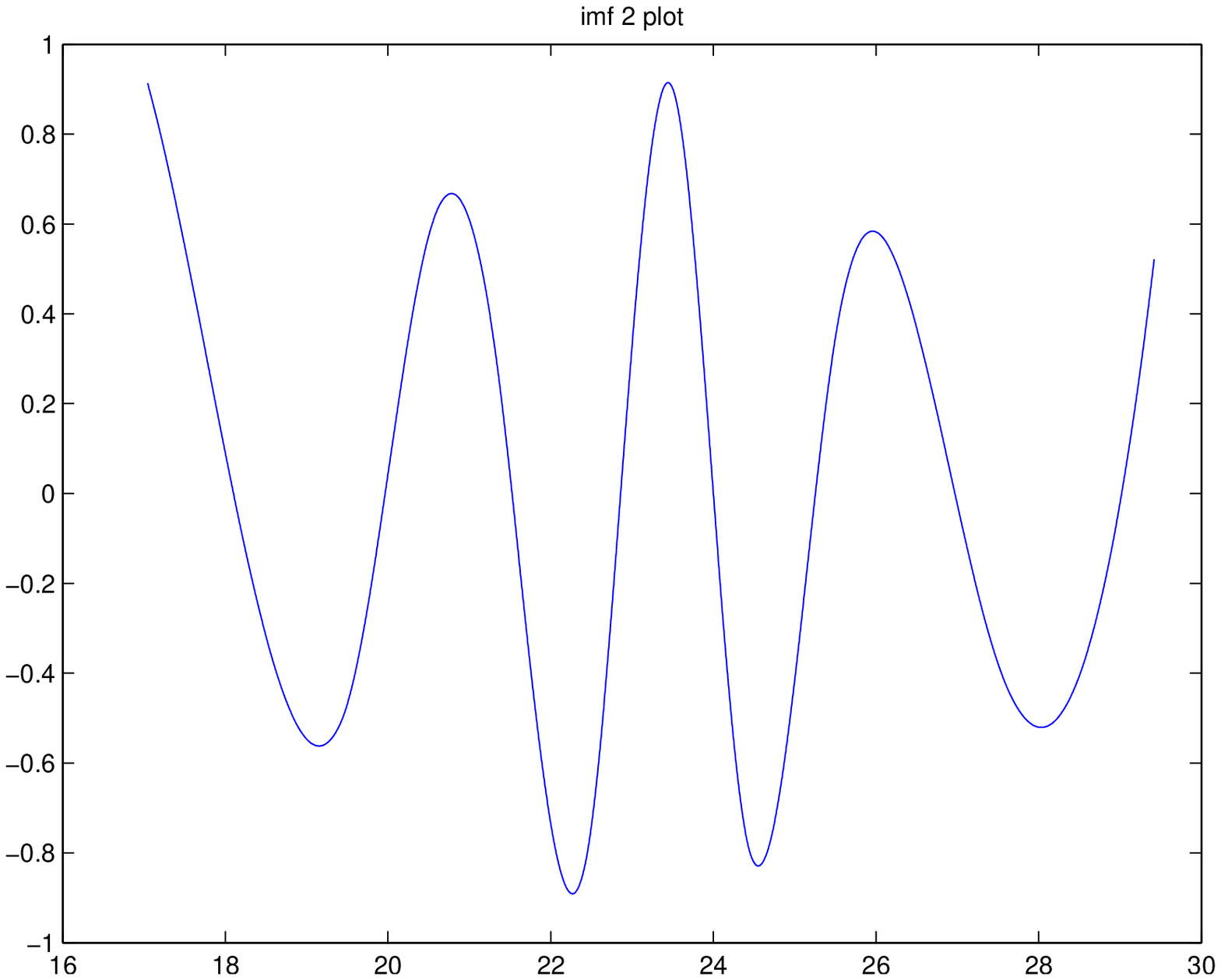}
\includegraphics[scale=1,height=160mm,angle=0,width=180mm]{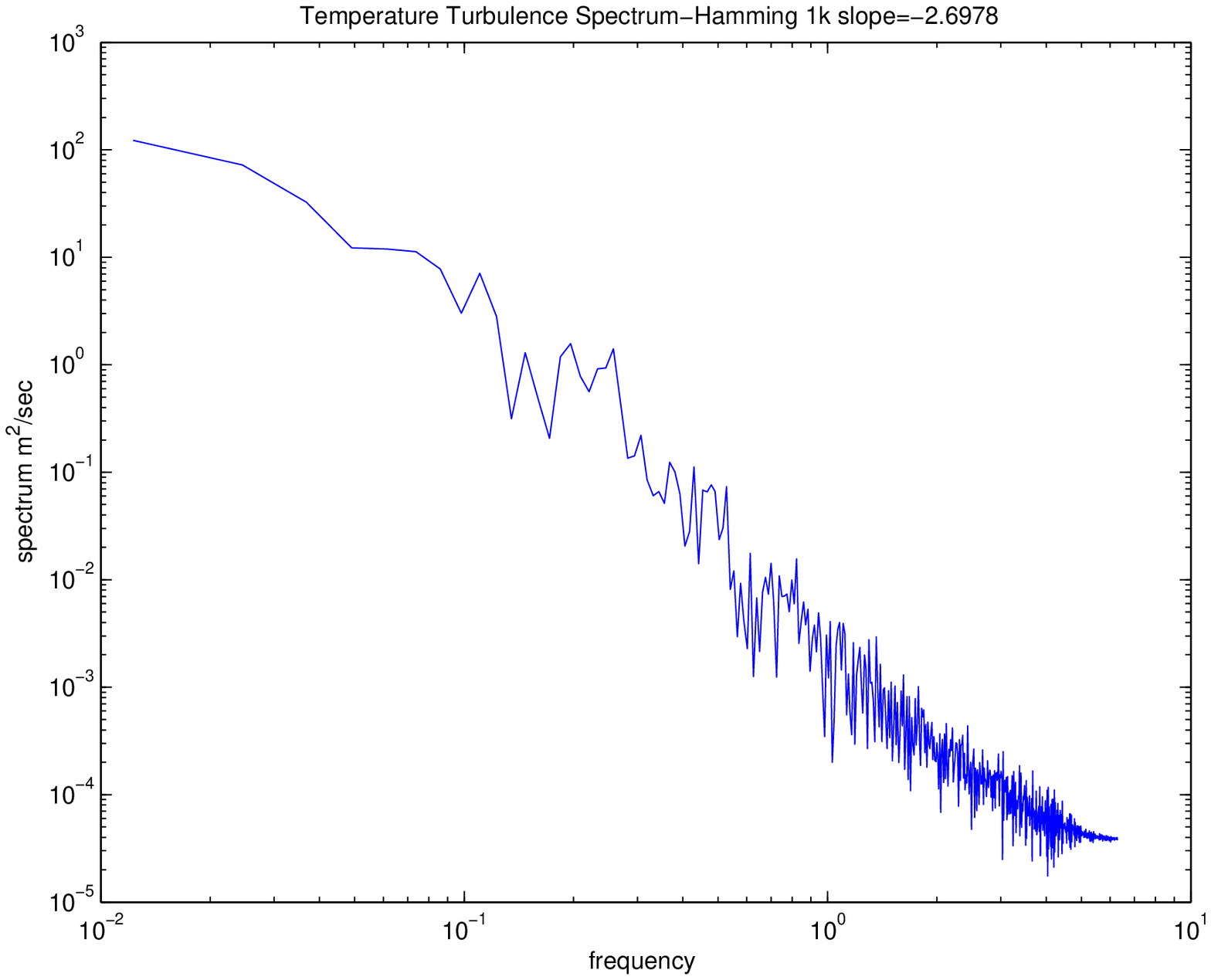}
\includegraphics[scale=1,height=160mm,angle=0,width=180mm]{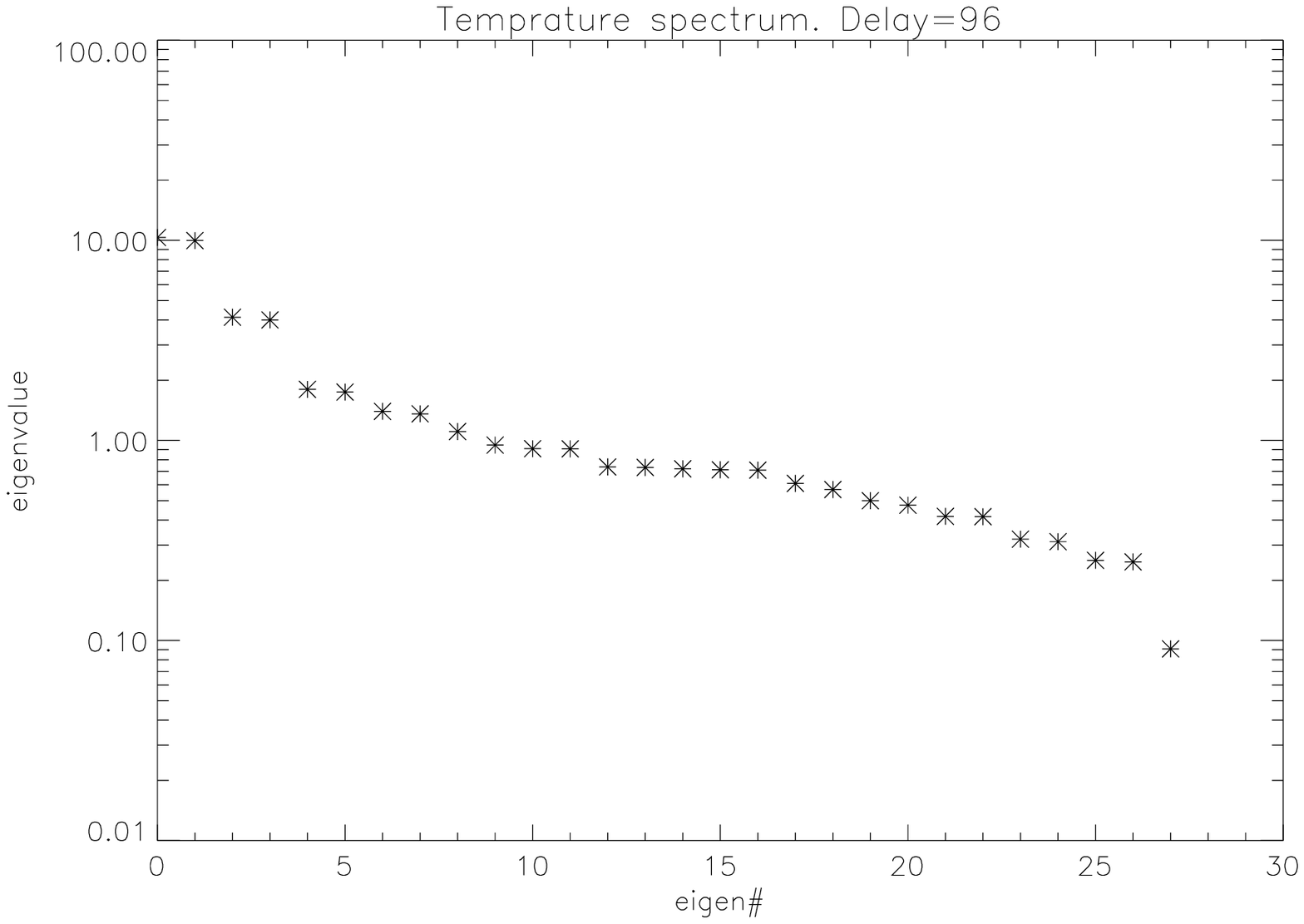}
\includegraphics[scale=1,height=160mm,angle=0,width=180mm]{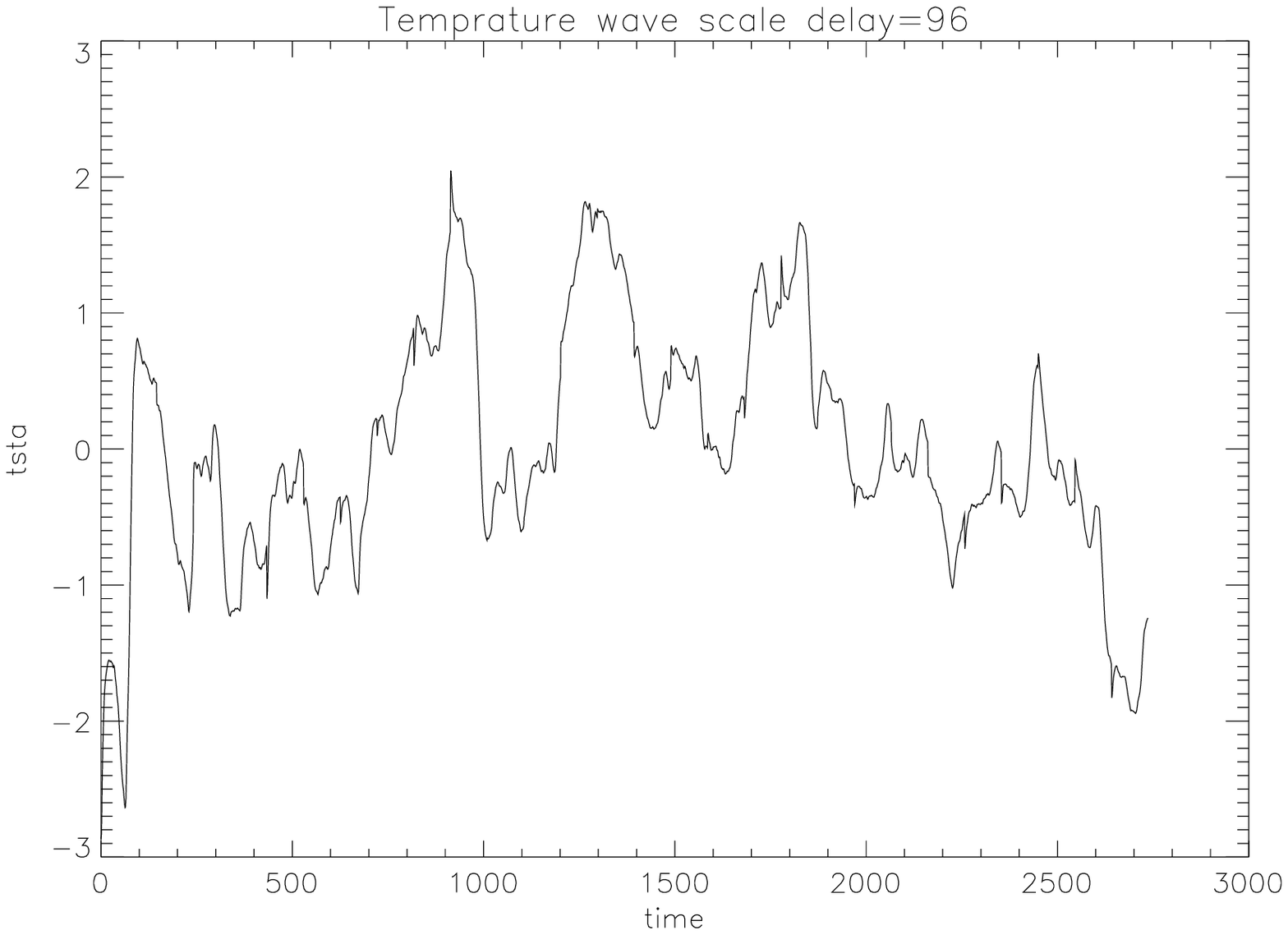}

\end{document}